\documentclass[12pt]{article}               % The use of LaTeX2e is preferred.
\usepackage{geometry}
\usepackage{cite}
\usepackage{multirow}
\usepackage{amsmath}
\usepackage{amssymb}
\usepackage[usenames]{color}
\usepackage{graphicx}          % Include this line if your
 \usepackage{amsthm}                              % document contains figures,
\usepackage{epsfig}    % or this line, depending on which
 \usepackage{helvet}	
 
\newtheorem{definition}{Definition}
\newtheorem{theorem}{Theorem}

\newtheorem{remark}{Remark}
\newtheorem{problem}{Problem}
\newtheorem{lemma}{Lemma}

\newtheorem{statement}{Statement}
\title{\large \bf
Extremum Seeking Approach for Nonholonomic Systems\\ with Multiple Time Scale Dynamics\thanks{This work was supported in part by the German Research Foundation (projects GR 5293/1-1 and ZU 359/2-1).\newline
$^1$Institute of Mathematics, Alpen-Adria University of Klagenfurt, Austria\newline
        {\tt\small viktoriia.grushkovskay@aau.at}
\newline
$^2$Max Planck Institute for Dynamics of Complex Technical Systems,  Magdeburg, Germany\newline
{\tt\small zuyev@mpi-magdeburg.mpg.de}
\newline
$^3$Otto von Guericke University Magdeburg, Germany
\newline
$^4$Institute of Applied Mathematics and Mechanics, National Academy of Sciences of Ukraine        
}}

\author{Victoria Grushkovskaya$^{1,4}$ \and Alexander Zuyev$^{2,3,4}$
}

\date{}
\begin{document}

\maketitle

\begin{abstract}
In this paper, a class of nonlinear driftless control-affine systems satisfying the bracket generating condition is considered.
A gradient-free optimization algorithm is developed for the minimization of a cost function along the trajectories of the controlled system.
The algorithm comprises an approximation scheme with fast oscillating controls for the nonholonomic dynamics and a model-free extremum seeking component with respect to the output measurements.
Exponential convergence of the trajectories to an arbitrary neighborhood of the optimal point is established under suitable assumptions on time scale parameters of the extended system.
The proposed algorithm is tested numerically with the Brockett integrator for different choices of generating functions.

\quad


\end{abstract}

%===============================================================================

\section{Introduction}

Extremum seeking theory aims at designing universal control algorithms which steer the trajectories of dynamical systems with uncertainties to the minimum (or maximum) of a cost function whose analytical representation may be partially or completely unknown.
The first results in this direction date back to the twenties of the last century, while the first thorough analysis of the stability properties of extremum seeking systems has been carried out only in the early 2000s, cf.~\cite{Kr00}.  Since then, many new extremum seeking algorithms and their applications have been developed (see, e.g., \cite{Krst03,Tan06,Nes10,Fu11,Liu12,Durr13,Har13,Guay15,Ben16,GE16,EMGG17,Pov17,SK17,SD17,GZE18,Gu18,Lab19,Mand19}). A special place in these extremum seeking studies is given to nonlinear systems with dynamic input-output maps of the form
\begin{equation}\label{sys_gen}
\begin{aligned}
&\dot x=f(x,\xi),\quad x\in\mathbb R^n, \xi\in\mathbb R^m,\,f:\mathbb R^n\times\mathbb R^m\to \mathbb R^n,\\
&y=h(x,\xi),\quad y\in\mathbb R^p,\,h:\mathbb R^n\times\mathbb R^m\to\mathbb R^p.
\end{aligned}
\end{equation}
The classical extremum seeking problem statement for system~\eqref{sys_gen} is to define the input $\xi$ in such a way that the output of system~\eqref{sys_gen} is optimized in the sense of minimization (or maximization) of an output-dependent cost function $J:\mathbb R^p\to\mathbb R$.
In this direction one can mention, e.g., the papers by~\cite{Kr00,Tan06,Ghaf12,Guay15,Har16,Durr17,Gu18}. Typically,
extremum seeking approaches for~\eqref{sys_gen} are based on the construction of a dynamic extension $\dot \xi =g\big(J(y),t\big)$, where $g:\mathbb R\times [0,\infty)\to \mathbb R^m$ is chosen to ensure the desired vicinity of the trajectories of~\eqref{sys_gen} to an optimal point. The analysis of the resulting system relies on singular perturbation theory and requires that system~\eqref{sys_gen} admits a steady-state $x=\ell(\xi)$, which is asymptotically stable for each fixed value of $\xi$. Furthermore, a crucial assumption in such studies is the existence of certain Lyapunov function for system~\eqref{sys_gen}.
However, there are many important classes of systems which do not admit a control Lyapunov function with desired properties.

In this paper, we consider a class of nonholonomic systems governed by driftless control-affine systems, in which the number of inputs can be significantly smaller than the number of state variables. In general, the linearization of these systems is not controllable. Moreover, as it was proved in the famous work by~\cite{Bro83}, such nonholonomic systems cannot be stabilized by a continuous feedback law. To stabilize such systems one can use, e.g., discontinuous~(e.g., \cite{ast94,Clar97}) or time-varying feedback laws~(e.g., \cite{ZuSIAM,GZ18}). Consequently, the resulting closed-loop system becomes discontinuous or non-autonomous and, in general, does not admit a regular Lyapunov function of the form $V(x)$.

The goal of our paper is to construct extremum seeking controls for a class of nonholonomic systems with time-varying inputs adapted from~\cite{GZ18}. We propose a novel solution of the extremum seeking problem for nonholonomic systems  based on combination of stabilizing strategies for nonholonomic systems and gradient-free extremum seeking controllers.
Although the main idea of our control design approach is inspired by singular perturbation techniques, we do not apply them directly in the proof. Instead, we propose a novel approach for dynamic stabilization of nonholonomic systems and generalize the techniques introduced in~\cite{GZE18} to systems with multiple time scales.

 The rest of this paper is organized as follows. In Section~\ref{sec_prel}, we introduce basic notations, formulate the problem statement, and describe the main idea of our control design approach. Section~\ref{sec_main} provides the main results of the paper, which are illustrated with an example in Section~\ref{sec_exmpl}. Section~\ref{sec_concl} contains  concluding remarks. Some auxiliary statements are given in Appendix~A, and the proof of the main result is contained in Appendix~B.
\section{Preliminaries}~\label{sec_prel}
\subsection{ Notations and Definitions}

%Throughout the paper, we exploit the following notations:
$\delta_{ij}$ is the Kronecker delta;
%$e_j$ -- unit vector with non-zero $j$-th entry;

${\rm dist}(x,S)$ is the Euclidian distance between an $x\in\mathbb R^{n}$ and an $S\subset\mathbb R^{n}$;

$B_\delta(x^*)$ is a $\delta$-neighborhood of an $x^*\in \mathbb R^n$;

$\partial M$, $\overline M$ is the boundary and the closure of a set $M\subset\mathbb R^n$, respectively; $\overline M= M\cup \partial M$;

$|S|$ is the cardinality of a set $S$;

$\mathcal K$ is  the class of  continuous strictly increasing functions $\varphi:\mathbb R^+\to\mathbb R^+$ such that $\varphi(0)=0$;

$[f,g](x)$ is the Lie bracket of vector fields $f,g:\mathbb R^n\to\mathbb R^n $ at a point $x\in\mathbb R^n$,   $[f,g](x)=L_fg(x)- L_gf(x)$, where  $ L_gf(x)=\lim_{s\to0}\dfrac{f(x+sg(x))-f(x)}{s}$.

Similarly to~\cite{Clar97,ZuSIAM}, we  exploit the sampling approach for the stabilization of nonholonomic systems.
Given an $\varepsilon{>}0$, we define the partition $\pi_\varepsilon$ of $[0,+\infty)$ into the intervals
$$
I_j=[t_j,t_{j+1}),\;t_j=\varepsilon j, \ j\in\mathbb N{\cup}\{0\}.
$$

\begin{definition}
Assume given a feedback $u=\varphi(x,\xi,t)$, $\varphi:D\times D\times[0,+\infty)\to\mathbb R^m$, $\varepsilon>0$, and $x^0,\xi^0\in D\subseteq \mathbb R^n$. \emph{A $\pi_\varepsilon$-solution} of the system
\begin{equation}
\begin{aligned}
&\dot x=f(x,u),\ \dot \xi =g(x,\xi,t),\, x,\xi\in D\subseteq \mathbb R^n, u\in\mathbb R^m,
\end{aligned}
\label{extended}
\end{equation}
 corresponding to $(x^0,\xi^0,\varphi)$, is an absolutely continuous function  $(x^\top(t),\xi^\top(t))^\top\in D\times D$, defined for $t\in[0,+\infty)$, which satisfies the initial conditions  $x(0)=x^0$, $\xi(0)=\xi^0$ and the differential equations
$$
\begin{aligned}
&\dot x(t)=f\big(x(t), \varphi(x(t_j),\xi(t_j),t)\big), \quad t\in I_j=[t_j,t_{j+1}),\\
&\dot \xi(t) = g(x(t),\xi(t),t)\text{ for each }j=0,1,2,\dots\ .
\end{aligned}
$$

\end{definition}
The above definition will be applied for the stabilization of nonholonomic systems using the approach of~\cite{ZuSIAM,GZ18}. However, the extremum seeking scheme proposed in this paper can also be used for output stabilization of systems with well-defined classical solutions.

\subsection{ Problem statement \& Main idea}
Consider a class of nonholonomic systems governed by driftless control-affine equations with single output:
\begin{equation}\label{sys1}
\begin{aligned}
&\dot x=\sum_{i=1}^mu_if_i(x),\\
& y=J(x),
\end{aligned}
\end{equation}
where $x=(x_1,\dots,x_n)^\top\in D{\subseteq} \mathbb R^n$ is the state, $x(0)=x^0\in D$, $u=(u_1,\dots,u_m)^\top\in\mathbb R^m$ is the control, $m<n$, $y\in\mathbb R$ is the output of the system, $J:D\to\mathbb R$ is the cost function, and the vector fields $f_i:D\to \mathbb R^n$ are linearly independent.
Let the following rank condition be satisfied in $D$:
\begin{equation}\label{rank}
\begin{aligned}
 {\rm span}\big\{f_{i}(x), [f_{j_1},f_{j_2}](x) \,|\,i\in S_1,(j_1,j_2)\in S_2\big\}=\mathbb{R}^n,
\end{aligned}
\end{equation}
where $S_1\subseteq \{1,2,...,m\}$ and  $S_2\subseteq \{1,2,...,m\}^2$ are some sets of indices, $|S_1|+|S_2|=n$.
We study the following extremum seeking problem:
\begin{problem}
Let $J\in C^2( D;\mathbb R)$ be a strongly convex function, and let $x^*\in D$ be such that $J(x)>J(x^*)$ for all $x\in D\setminus\{x^*\}$. The goal is to construct a control law $u=u(t,x,J(x))$ such that the trajectories $x(t)$ of system~\eqref{sys1} with the initial conditions from $D$ tend asymptotically to an arbitrary small neighborhood of $x^*$.
\end{problem}
The main idea of the control algorithm proposed in this paper can be described in two stages:

(1) \underline{Model-based stabilizing component.}
  For each value $\xi\in D$, we construct time-periodic fast oscillating control laws with \emph{state-dependent} coefficients to ensure that the corresponding steady-state $x=\xi$ of~\eqref{sys1} is asymptotically (and even exponentially) stable.
  Further we assume that $\xi(t)$ evolves according to certain differential equations, so the result of~\cite{ZuSIAM,GZ18} cannot be directly applied for establishing stability properties of the extended system~\eqref{extended}.  Note that, in general, \eqref{sys1} does not admit a control Lyapunov function. Instead, we will prove that with the proposed choice of the control $u$ the trajectory $x(t)$ remains in a sufficiently small neighborhood of $\xi(t)$ for $t\in[0,\infty)$.  These controls are model-based, i.e. the dynamics (control vector fields) and the coordinates of the system are assumed to be known, but not the analytical expression of $J$ and the optimal point $x^*$. We will apply sampling controllers, that is the solutions of~\eqref{sys1} will be defined in the sense of Definition~1.
  
(2) \underline{Model-free extremum seeking component.} To optimize the state $x=\xi$ with respect to minimizing the cost function $J(x)$ along the trajectories of~\eqref{sys1}, we construct a dynamic extension $\dot \xi=g(y,t)$, where $g(y,t)$ is taken in the form of fast oscillating time-periodic functions  with \emph{output-dependent} coefficients from~(\cite{GZE18}).  Thus, this part of the controller is model-free.
      %Although in this paper we assume that the existence of the steady-state $x=\xi$ is ensured by the controls from (\cite{ZuSIAM,GZ18}), the proposed extremum-seeking control design can be applied to other classes of systems under the assumption that $\|x(t)-\xi(t)\|$ is small enough for $t\in[0,\infty)$.
\begin{remark}
{In Problem~1, we assume that the cost function $J$ depends only on the state variable $x$, but not on the control input $u$. This assumption is not crucial and is made in order to simplify the proof. Besides, if $J$ depends only on   $u$, the stability properties directly follow from~(\cite{GZE18}) and~(\cite{GZ18}) with the same proof techniques.
}\end{remark}
\normalsize
\vspace{-2em}
\section{Main results}~\label{sec_main}
\subsection{Control design}
In this section, we formalize the control algorithm announced in Subsection~2.2. Namely, the overall system has the following form:
\begin{subequations}
\begin{align}
&\dot x=\sum_{i=1}^mu_if_i(x),\, x(0)=x^0\nonumber,\\
& u_i=\varphi_i^\varepsilon(x,\xi,t),\label{nonhA}\\
&y=J(x),\nonumber\\
&\dot \xi=g(y,t),\, g(y,t)=\sum_{j=1}^{2n}g_{j}(y)v_{j}^\mu(t)e_j,\;\quad\;\,\xi(0)=x^0.\label{nonhB}
\end{align}\label{nonh}
\end{subequations}
In~\eqref{nonhA}, the stabilizing component $u_i=\varphi_i^\varepsilon(x,\xi,t)$ is
{\begin{align}
&\varphi_i^\varepsilon(x,\xi,t)=\sum_{i_1\in S_1}a_{i_1}(x,\xi)\delta_{i i_1}\label{cont}\\
  &+\sqrt{\dfrac{4\pi}{\varepsilon}}\sum_{(i_1,i_2)\in S_2}{\sqrt{\kappa_{i_1i_2}|a_{i_1i_2}(x,\xi)|}} \Big(\delta_{ii_1}{\rm sign}(a_{i_1,i_2}(x,\xi))\cos{\dfrac{2\pi \kappa_{i_1i_2}}{\varepsilon}}t+\delta_{ii_2}\sin{\dfrac{2\pi \kappa_{i_1i_2}}{\varepsilon}}t\Big).\nonumber
\end{align}
}
 Here $\kappa_{i_1i_2}\in\mathbb N$, $\kappa_{i_1i_2}\ne\kappa_{i_3i_4}$ for all $(i_1,i_2)\ne(i_3,i_4)$,    and 
 $$a(x,\xi)=\Big((a_{i_1}(x,\xi))_{i_1\in S_1}\  ( a_{i_1i_2}(x,\xi))_{(i_1,i_2)\in S_2}\Big)^\top\in\mathbb R^n$$
 is defined as
{\begin{equation}\label{a}
a(x,\xi)=-  \gamma_1 \mathcal F^{-1}(x) (x-\xi)
\end{equation}}
with $\mathcal F^{-1}(x)$ being the $n\times n$ matrix inverse to
  $$
  \mathcal F(x)= \Big(\big(f_{j_1}(x)\big)_{j_1\in S_1}\ \ \big([f_{j_1},f_{j_2}](x)\big)_{(j_1,j_2)\in S_2}\Big),
  $$
and the control gain $\gamma_1{>}0$ to be defined later in the proof of the main result. 

Such a choice of  $u_i$ is aimed to ensure that the trajectories $x(t)$ are  close enough to $\xi(t)$ for all $t\ge 0$ and all initial conditions $x(0)$. Note that the rank condition~\eqref{rank} implies nonsingularity of $\mathcal F(x)$ for any $x\in D$.

 In~\eqref{nonhB}, $g(y,t)$ is the extremum seeking component. Here  $e_j$ denotes the unit vector in $\mathbb R^n$ with non-zero $j$-th entry if $j\le n$, and non-zero $(j-n)$-th entry if $n+1\le j\le 2n$, the functions $g_j,g_{j+n}$ have to satisfy the relation
 $$
 [g_j(z),g_{j+n}(z)]=-\gamma_2,\;\gamma_2>0, \ j=\overline{1,n}.
 $$
 For example, the choice $g_{j+n}(z)=-\gamma_2 g_{j}(z)\int{\dfrac{dz}{g_{j}(z)^2}}$ was proposed in~\cite{GZE18}. In this paper, we propose to parameterize  the functions $g_j,g_{j+n}$ as
 \begin{equation}\label{class}
  \begin{aligned}
 &g_j(z)=r_j(z)\sin\phi_j(z),\,g_{j+n}(z)=r_j(z)\cos\phi_j(z),\\
 &\text{with }r_j,\phi_j\text{ such that }r_j^2(z)\phi'_j(z)\equiv \gamma_2.
 \end{aligned}
 \end{equation}
The discrete-time version of the above parametrization has also been used by~\cite{FG19}.

 %according to~\cite{GZE18},
%\begin{equation}\label{class}
%g_{j+n}(z)=-\gamma_2 g_{j}(z)\int{\dfrac{dz}{g_{j}(z)^2}},\;\gamma_2>0, \ j=\overline{1,n},
%\end{equation}
Next, the inputs $v_{j}^\mu(t)$ are given by
 $$
  v_{j}^\mu(t)=\left\{
\begin{aligned}
\sqrt{\dfrac{4\pi  k_j}{\mu}}\cos\dfrac{2\pi k_jt}{\mu} &\text{ for }j=\overline{1,n},\\
 \sqrt{\dfrac{4\pi  k_{j-n}}{\mu}}\sin\dfrac{2\pi k_{j-n}t}{\mu} &\text{ for }j=\overline{n+1,2n},
\end{aligned}
\right.
 $$
where $\mu>0$,  $k_j\in\mathbb N$, $k_{j_1}\ne k_{j_2}$ for all $j_1\ne j_2$.

\begin{remark}
  {Although the choice of $g_j,g_{j+n}$ in~\eqref{class} may look rather artificial, there are many extremum seeking systems whose control vector fields satisfy this relation. For example, the functions $g_{j}(z)=z$, $g_{j+n}(z)=1$ have been exploited by~\cite{Durr13a,Durr17};  $g_{j}(z)=\sin\ z$, $g_{j+n}(z)=\cos\ z$ by~\cite{SK17};  $g_{j}(z)=\sqrt z\sin(\ln z)$, $g_{j+n}(z)=\sqrt z\cos(\ln z)$ by~\cite{SD17}; $g_j(z)=\sqrt{\dfrac{1-e^{-z}}{1+e^{z}}}\sin(e^{z}+2\ln(e^{z}-1))$, $g_{j+n}(z)=\sqrt{\dfrac{1-e^{-z}}{1+e^{z}}}\cos(e^{z}+2\ln(e^{z}-1))$ by~\cite{GZE18}. One more example will be given in Section~\ref{sec_exmpl}.}
\end{remark}

\subsection{ Stability conditions}
Assume that the cost  function $J\in C^2(D;\mathbb R)$ satisfies the following properties in $D$:
{  \begin{align}
\sigma_{11}\|x-x^*\|^2\le J(x)-J^* &\le\sigma_{12}\|x-x^*\|^2,\nonumber\\
\sigma_{21}(J(x)-J^*)\le \|\nabla J(x)\|^2&\le\sigma_{22}(J(x)-J^*),\label{J}\\
&\Big\|\dfrac{\partial^2J(x)}{\partial x^2}\Big\|\le\sigma_3,\nonumber
  \end{align}}
with $x^*\in D$ and some positive constants $\sigma_{11}$, $\sigma_{12}$, $\sigma_{21}$, $\sigma_{22}$, $\sigma_3$.
The main result of this paper is as follows.
\begin{theorem}~\label{thm_step1}
{Given system~\eqref{sys1} and a function $J\in C^2(D;\mathbb R)$ satisfying~\eqref{J}, assume that:
\begin{itemize}
\item the  vector fields $f_i\in C^2(D;\mathbb R^n)$ in~\eqref{sys1}  satisfy~\eqref{rank} in  $D$, and
 there is an  $\alpha{>}0$ such that $\|\mathcal F^{-1}(x)\|\le \alpha \text{ for all }x\in D;$\\
\item $g_{j}(J(\cdot))\in C^2(D\setminus\{x^*\};\mathbb R)$, $L_{g_{j}}g_{i}(J(\cdot))$, $L_{g_{l}}L_{g_{j}}g_{i}(J(\cdot))\in C(D;\mathbb R)$ for all  $i,j,l=\overline{1,2n}$;\\
\item for any compact  $D'\subseteq D$, there are  $L_g, L_{2g}, M_{3g}\ge 0$ s.t.
$$
 \begin{aligned}
&\|g_{i}\big(J(x)\big)-g_{i}\big(J(\xi)\big)\|\le L_g\|x-\xi\|,\\
&\|L_{\big(g_{j_2}(J(x))-g_{j_2}(J(\xi))\big)}g_{j_1}\big(J(\xi)\big)\|\le L_{2g}\|x-\xi\|,\\
&\|L_{g_{j_3}(J(x))} L_{g_{j_2}(J(\xi))}g_{j_1}\big(J(\xi)\big)\|\le M_{3g},\quad x,\xi\in D',\,i,j,l=\overline{1,2n}.
\end{aligned}
$$
\end{itemize}
 Then, for any $\delta\in\Big(0,\sqrt{{\sigma_{11}}/{\sigma_{12}}}{\rm dist}(x^*,\partial D)\Big)$ and any $\rho{>}0$,  there exist  $\bar\mu>0,\bar\gamma_1(\mu)>0$, and $\bar\varepsilon(\gamma_1,\mu)>0$  such that, for any $\mu\in(0,\bar\mu]$, $\gamma_1\in[\bar \gamma_1(\mu),\infty)$, and any $\varepsilon\in(0,\bar\varepsilon(\gamma_1,\mu)]$, each $\pi_\varepsilon$ solution of~\eqref{nonh} with  $u_i=\varphi_i^\varepsilon(x,\xi,t)$ defined by~\eqref{cont} and the initial conditions from $\overline{B_\delta(x^*)}$ satisfies
\begin{equation}\label{decay}
\|x(t)-x^*\|\le \beta\|x^0-x^*\|e^{-\lambda t}+\rho\text{ for all }t\in[0,\infty),
\end{equation}
with some $\beta,\lambda>0$.}
  \end{theorem}
The proof of this theorem is given in Appendix~B.
{
\begin{remark}
{  The proof of Theorem~\ref{thm_step1} represents a constructive procedure for choosing $\bar\mu$, $\bar\gamma_1(\mu)$, $\bar\varepsilon(\gamma_1,\mu)$, and clarifies the relation between these parameters and the coefficients $\beta$ and $\lambda$. We would like to underline that the proposed bounds are quite conservative. The crucial assumption is $\varepsilon<\mu$, which means that subsystem~\eqref{nonhA} oscillates faster than subsystem~\eqref{nonhB}. To simplify the proof, we also suppose that $\dfrac{\mu}{\varepsilon}\in\mathbb N$ and $x(0)=\xi(0)$, however the assertion of Theorem~1 can also be obtained without these assumptions.}
\end{remark}}
In order to have $\gamma_1$ independent on $\mu$, one may introduce an additional parameter $\eta$ which will ensure a ``slow'' dynamics of~\eqref{nonhB} (similarly to, e.g.,~\cite{Durr17}). This, however, will result in a slower convergence rate of the overall system to the optimal point. Namely, by taking $\tilde v_j^\mu(t):=\dfrac{1}{\eta}v_j^\mu\Big(\dfrac{t}{\eta}\Big)$ in~\eqref{nonhB} and keeping the conditions of Theorem~\ref{thm_step1}, one can prove the following statement: 

\emph{ For any $\delta\in\Big(0,\sqrt{\dfrac{\sigma_{11}}{\sigma_{12}}}{\rm dist}(x^*,\partial D)\Big)$ and any $\rho>0$,  there exist  $\bar\mu>0$, $\bar\varepsilon(\mu)>0$, and $\bar\eta(\varepsilon,\mu)>0$  such that, for any $\mu\in(0,\bar\mu]$, $\varepsilon\in(0,\bar\varepsilon(\mu)]$ and $\eta\in[\bar \eta(\varepsilon,\mu),\infty)$, each $\pi_\varepsilon$-solution of~\eqref{nonh} with  $u_i=\varphi_i^\varepsilon(x,\xi,t)$ defined by~\eqref{cont} and the initial conditions from $\overline{B_\delta(x^*)}$ satisfies
$
\|x(t)-x^*\|\le \beta\|x^0-x^*\|e^{-\dfrac{\lambda t}{\eta}}+\rho\text{ for all }t\in[0,\infty),\,\beta,\lambda>0.
$}

Similarly to~\cite{GZE18}, the behavior of the solutions of~\eqref{sys1} can be improved by generating  $g_j$ vanishing at $x^*$. We will illustrate this feature with an example in the next section.

\section{Example}~\label{sec_exmpl}
\begin{figure*}[ht]
 \begin{minipage}{0.49\linewidth}
\begin{center}
\includegraphics[width=1\linewidth]{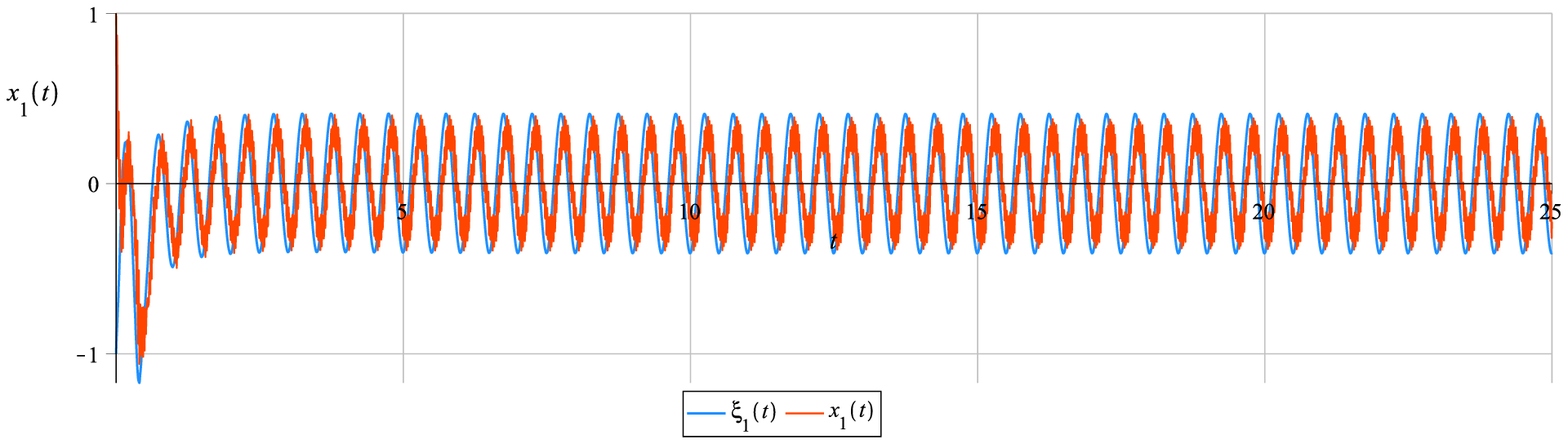}
\includegraphics[width=1\linewidth]{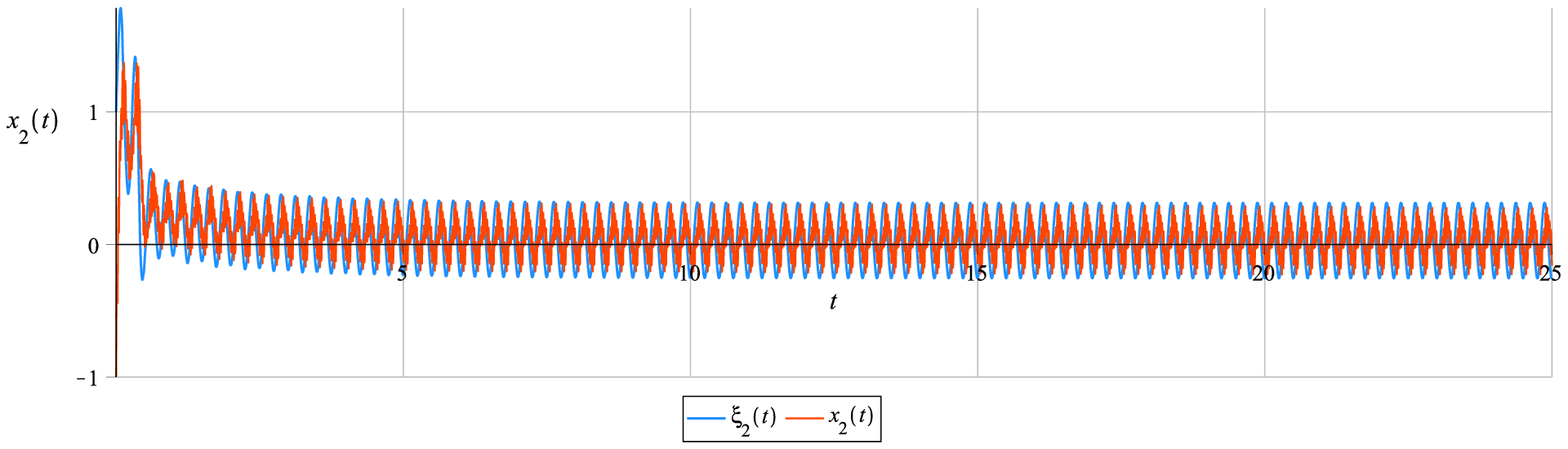}
\includegraphics[width=1\linewidth]{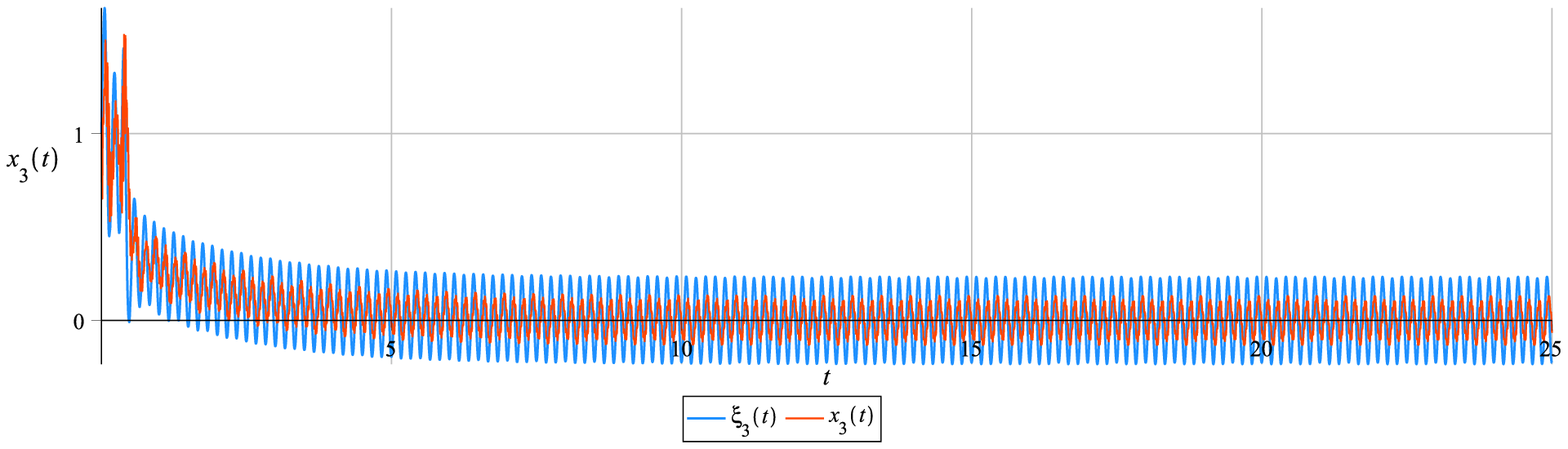}
\end{center}
 \end{minipage}\hfill
 \begin{minipage}{0.49\linewidth}
\begin{center}
\includegraphics[width=1\linewidth]{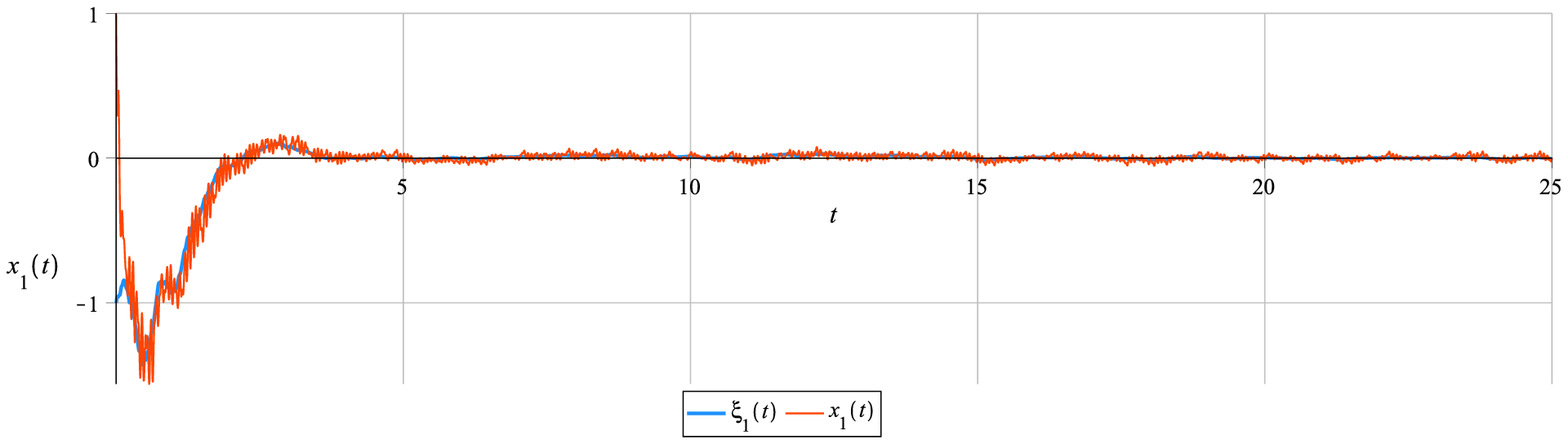}
\includegraphics[width=1\linewidth]{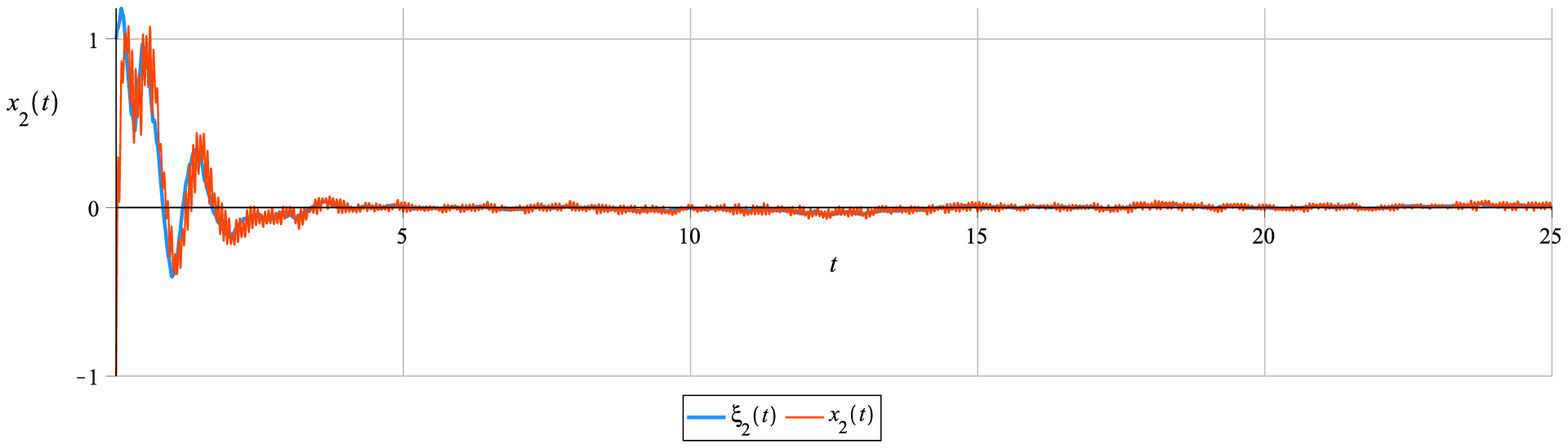}
\includegraphics[width=1\linewidth]{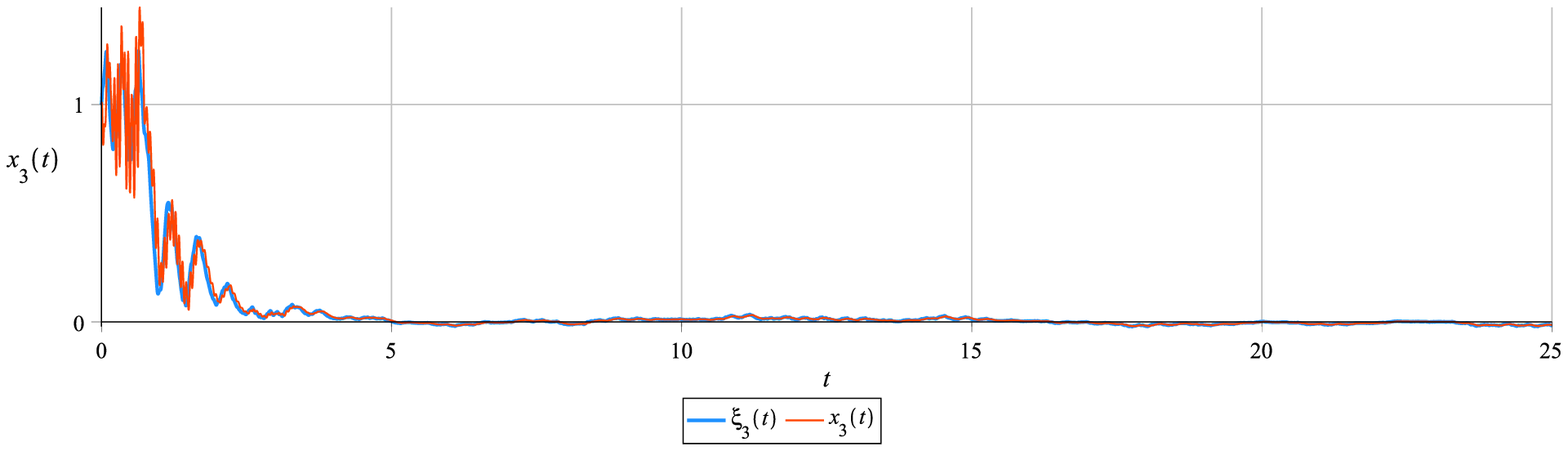}
\end{center}
 \end{minipage}
\caption{ Time-plots of the trajectories of system~\eqref{bro_x}--\eqref{bro_xi} with the generating functions~\eqref{class_durr} (left) and~\eqref{class_vv} (right).}\label{fig_bro}
\end{figure*}

As an example, consider the well-known Brockett integrator~(\cite{Bro83}):
\begin{equation}\label{bro_x}
\begin{aligned}
&  \dot x_1 =u_1,\ \dot x_2=u_2,\ \dot x_3=u_1x_2-u_2x_1.
\end{aligned}
\end{equation}
It is easy to see that, for all $x\in\mathbb R^3$, the vector fields $f_1=(1,0,x_2)^\top$ and $f_2=(0,1,-x_1)^\top$ of system~\eqref{bro_x} satisfy the rank condition~\eqref{rank} with $S_1=\{1,2\}$, $S_2=\{(1,2)\}$:
$
 {\rm span}\big\{f_1(x),\,f_2(x),\, [f_1,f_2](x)\big\}=\mathbb{R}^3\text{ for all }x\in\mathbb R^3;
$
thus, we may apply the control algorithm proposed in Section~3.1. Namely, we take
{
\begin{align}
&u_1=a_1(x,\xi)+\sqrt{{4\pi\kappa_{12}\big|a_{12}(x)\big|}/{\varepsilon}}\ {\rm sign}\big({a_{12}(x,\xi)}\big)\cos({2\pi\kappa_{12}t}/{\varepsilon}),\nonumber\\
&u_2=a_2(x,\xi)+\sqrt{{4\pi\kappa_{12}\big|a_{12}(x)\big|}/{\varepsilon}}\ \sin({2\pi\kappa_{12}t}/{\varepsilon}),\label{bro_u}\\
&a(x,\xi)=\big(a_1(x,\xi),a_2(x,\xi),a_{12}(x,\xi)\big)^\top=-  \gamma_1 \mathcal F^{-1}(x) (x-\xi)\nonumber\\
&\quad = -\gamma_1\left(
 x_1-\xi_1,
  x_2-\xi_2,
  \dfrac{1}{2}\big(-x_2\xi_1+x_1\xi_2-x_3+\xi_3\big)^\top
\right), \nonumber\\
&\dot \xi_j=\sqrt{{4\pi k_j}/{\mu}}\Big(g_1(y)\cos({2\pi k_jt}/{\mu})+g_2(y)\sin({2\pi k_jt}/{\mu})\Big)e_j,\label{bro_xi}
\end{align}
 $j=1,2,3$}. In this example, we take $y=J(x)=\|x\|^2$, $\gamma_1=20$, $\gamma_2=1$, $\kappa_{12}=4$, $k_1=1$, $k_2=2$, $k_3=3$, and consider two types of functions $g_1,g_2$. The results of numerical simulations with the functions from~\cite{Durr17},
\begin{equation}\label{class_durr}
g_1(z)=z,\,g_2(z)=1,
\end{equation}
are depicted on Fig.~\ref{fig_bro} (left). Here $\varepsilon=0.1$ and $\mu=0.5$.

To improve the qualitative behavior of~\eqref{bro_x}--\eqref{bro_xi}, we can apply another pair of  the generating functions satisfying~\eqref{class}, which vanish when $J$ takes its minimal value, e.g.,
{
\begin{align}
 & g_1(z)=\sqrt{\tanh{z}/{2}}\sin\big(2\ln(e^z-1)-z\big),\label{class_vv}\\
 &  g_2(z)=\sqrt{\tanh{z}/{2}}\cos\big(2\ln(e^z-1)-z\big)\text{ if }z>0,g_1(0)=g_2(0)=0.\nonumber
\end{align}
}
\vspace{-1em}

In this case, we took $\varepsilon=0.25$, $\mu=1$. Note that, unlike the results of~\cite{GZE18}, the trajectories of~\eqref{bro_x}--\eqref{bro_xi} exhibit non-vanishing oscillations in a neighborhood of the extremum point (which are, however, considerably smaller than with the functions~\eqref{class_durr}) (see Fig.~\ref{fig_bro}, right). Thus, an interesting question is whether it is
possible to achieve asymptotic stability in the sense of Lyapunov with the proposed control algorithm.

In both case, we take the initial conditions $x(0)=(1,-1,1)^\top$, $\xi(0)=(-1,1,1)^\top$ to illustrate that the proposed approach can be applied also for $x^0\ne \xi^0$.

\section{Conclusions \& Future work}~\label{sec_concl}
To simplify the presentation, we consider only the class of nonholonomic systems~\eqref{sys1} satisfying one-step bracket generating condition in this paper, i.e. we assume that the vector fields together with their Lie brackets span the whole $n$-dimensional space at each state $x\in D\subseteq {\mathbb R}^n$.
Another hypothesis is put in~\eqref{J}, so that the cost $J$ possesses properties of a quadratic function.  This hypothesis is introduced in order not to overcomplicate the proof of the main results.
It should be emphasized that information about the analytical expression of $J$ and its minimizer $x^*$ is not required for the control design.
Furthermore, all the constants  in~\eqref{J} may also be unknown.
In future work, we expect to address broader classes of cost functions possessing polynomial convergence properties, similarly to the results of~\cite{GZE18}.
We also plan to extend the proposed control design approach to nonholonomic systems under higher order controllability conditions with iterated Lie brackets.

\bibliographystyle{plain}
\bibliography{biblio_es}             % bib file to produce the bibliography

\begin{thebibliography}{10}

\bibitem{ast94}
A.~Astolfi.
\newblock On the stabilization of nonholonomic systems.
\newblock In {\em Proc. 33rd IEEE Conference on Decision and Control},
  volume~4, pages 3481--3486, 1994.

\bibitem{Ben16}
M.~Benosman.
\newblock {\em Learning-Based Adaptive Control: An Extremum Seeking
  Approach--Theory and Applications}.
\newblock Butterworth-Heinemann, 2016.

\bibitem{Bro83}
R.~W. Brockett.
\newblock Asymptotic stability and feedback stabilization.
\newblock {\em Differential Geometric Control Theory}, pages 181--191, 1983.

\bibitem{Clar97}
F.~H. Clarke, Y.~S. Ledyaev, E.~D. Sontag, and A.~I. Subbotin.
\newblock Asymptotic controllability implies feedback stabilization.
\newblock {\em IEEE Tran on Automatic Control}, 42(10):1394--1407, 1997.

\bibitem{Durr17}
H.-B. D{\"u}rr, M.~Krsti{\'c}, A.~Scheinker, and C.~Ebenbauer.
\newblock Extremum seeking for dynamic maps using {L}ie brackets and singular
  perturbations.
\newblock {\em Automatica}, 83:91--99, 2017.

\bibitem{Durr13a}
H.-B. D\"{u}rr, M.~S. Stankovi\'{c}, C.~Ebenbauer, and K.H. Johansson.
\newblock Lie bracket approximation of extremum seeking systems.
\newblock {\em Automatica}, 49:1538--1552, 2013.

\bibitem{Durr13}
H.~B. D\"{u}rr, C.~Zeng, and C.~Ebenbauer.
\newblock Saddle point seeking for convex optimization problems.
\newblock {\em Proc. 9th IFAC Symposium on Nonlinear Control Systems}, pages
  540--545, 2013.

\bibitem{EMGG17}
C.~Ebenbauer, S.~Michalowsky, V.~Grushkovskaya, and B.~Gharesifard.
\newblock Distributed optimization over directed graphs with the help of {L}ie
  brackets.
\newblock {\em IFAC-PapersOnLine}, 50(1):15343--15348, 2017.

\bibitem{FG19}
J.~Feiling, C.~Labar, V.~Grushkovskaya, E.~Garone, M.~Kinnaert, and
  C.~Ebenbauer.
\newblock Extremum seeking algorithms based on non-commutative maps.
\newblock {\em IFAC-PapersOnLine}, 52(16):688--693, 2019.

\bibitem{Fu11}
L.~Fu and {\"U}.~{\"O}zg{\"u}ner.
\newblock Extremum seeking with sliding mode gradient estimation and asymptotic
  regulation for a class of nonlinear systems.
\newblock {\em Automatica}, 47(12):2595--2603, 2011.

\bibitem{Ghaf12}
A.~Ghaffari, M.~Krsti{\'c}, and D.~Ne{\v{s}}I{\'c}.
\newblock Multivariable newton-based extremum seeking.
\newblock {\em Automatica}, 48(8):1759--1767, 2012.

\bibitem{GE16}
V.~Grushkovskaya and C.~Ebenbauer.
\newblock Multi-agent coordination with {L}agrangian measurements.
\newblock {\em IFAC-PapersOnLine}, 49(22):115--120, 2016.

\bibitem{GZ18}
V.~Grushkovskaya and A.~Zuyev.
\newblock Obstacle avoidance problem for second degree nonholonomic systems.
\newblock In {\em Proc. 57th IEEE Conf. on Decision and Control}, pages
  1500--1505, 2018.

\bibitem{GZE18}
V.~Grushkovskaya, A.~Zuyev, and C.~Ebenbauer.
\newblock On a class of generating vector fields for the extremum seeking
  problem: Lie bracket approximation and stability properties.
\newblock {\em Automatica}, 94:151--160, 2018.

\bibitem{Gu18}
M.~Guay and K.~T. Atta.
\newblock Dual mode extremum-seeking control via {L}ie-bracket averaging
  approximations.
\newblock In {\em Proc. 2018 Annual American Control Conference}, pages
  2972--2977, 2018.

\bibitem{Guay15}
M.~Guay and D.~Dochain.
\newblock A time-varying extremum-seeking control approach.
\newblock {\em Automatica}, 51:356--363, 2015.

\bibitem{Har13}
M.~Haring, N.~van~de Wouw, and D.~Ne{\v{s}}i{\'c}.
\newblock Extremum-seeking control for nonlinear systems with periodic
  steady-state outputs.
\newblock {\em Automatica}, 49(6):1883--1891, 2013.

\bibitem{Har16}
M.~A. Haring and T.~A. Johansen.
\newblock Asymptotic stability of perturbation-based extremum-seeking control
  for nonlinear systems.
\newblock {\em IEEE Trans. Autom. Control}, 62(5):2302--2317, 2017.

\bibitem{Krst03}
M.~Krsti\'{c} and K.~B. Ariyur.
\newblock {\em Real-Time optimization by Extremum Seeking Control}.
\newblock Wiley-Interscience, 2003.

\bibitem{Kr00}
M.~Krsti{\'c} and H.-H. Wang.
\newblock Stability of extremum seeking feedback for general nonlinear dynamic
  systems.
\newblock {\em Automatica}, 36(4):595--601, 2000.

\bibitem{Lab19}
C.~Labar, E.~Garone, M.~Kinnaert, and C.~Ebenbauer.
\newblock Newton-based extremum seeking: A second-order {L}ie bracket
  approximation approach.
\newblock {\em Automatica}, 105:356--367, 2019.

\bibitem{Liu12}
S.-J. Liu and M.~Krsti\'{c}.
\newblock {\em Stochastic averaging and stochastic extremum seeking}.
\newblock Springer Science \& Business Media, 2012.

\bibitem{Mand19}
F.~Mandi{\'c}, N.~Mi{\v{s}}kovi{\'c}, and I.~Lon{\v{c}}ar.
\newblock Underwater acoustic source seeking using time-difference-of-arrival
  measurements.
\newblock {\em IEEE Journal of Oceanic Engineering}, 2019.

\bibitem{Nes10}
D.~Ne{\v{s}}i{\'c}, Y.~Tan, W.~H. Moase, and C.~Manzie.
\newblock A unifying approach to extremum seeking: Adaptive schemes based on
  estimation of derivatives.
\newblock In {\em Prc. 49th IEEE Conf. on Decision and Control}, pages
  4625--4630, 2010.

\bibitem{Pov17}
J.~I. Poveda and A.~R. Teel.
\newblock A framework for a class of hybrid extremum seeking controllers with
  dynamic inclusions.
\newblock {\em Automatica}, 76:113--126, 2017.

\bibitem{SK17}
A.~Scheinker and M.~Krsti{\'c}.
\newblock {\em Model-free stabilization by extremum seeking}.
\newblock Springer, 2017.

\bibitem{SD17}
R.~Suttner and S.~Dashkovskiy.
\newblock Exponential stability for extremum seeking control systems.
\newblock {\em IFAC-PapersOnLine}, 50(1):15464--15470, 2017.

\bibitem{Tan06}
Y.~Tan, D.~Ne{\v{s}}i{\'c}, and I.~Mareels.
\newblock On non-local stability properties of extremum seeking control.
\newblock {\em Automatica}, 42(6):889--903, 2006.

\bibitem{ZuSIAM}
A.~Zuyev.
\newblock Exponential stabilization of nonholonomic systems by means of
  oscillating controls.
\newblock {\em SIAM J. on Control and Optimization}, 54(3):1678--1696, 2016.

\bibitem{ZG17}
A.~Zuyev and V.~Grushkovskaya.
\newblock Motion planning for control-affine systems satisfying low-order
  controllability conditions.
\newblock {\em International Journal of Control}, 90(11):2517--2537, 2017.

\end{thebibliography}
                                                    % with bibtex (preferred)

\appendix
\section{Auxiliary results}
This section contains several technical results which be used for the proof of Theorem~\ref{thm_step1}.
\begin{lemma}~\label{lemma_x}
{ Let $D{\subseteq}\mathbb R^n$, $\xi(t){\in }D$, $t\in[0,\tau]$, be a  solution of the system
$$\dot \xi=\sum_{i=1}^l h_i(\xi)w_i(t),
$$
 and let the vector fields $h_i$ be Lipschitz continuous in $D$ with the Lipschitz constant $L$. Then
$$
      \|\xi(t)-\xi(0)\|\le t\nu\max_{1\le i\le l}\|h_i(\xi(0))\|e^{\nu Lt},\;t\in[0,\tau],
$$
    with $\nu=\max\limits_{t\in[0,\tau]}\sum_{i=1}^{l}|w_{i}(t)|$.}
  \end{lemma}
Lemma~\ref{lemma_x} follows from the  Gr\"{o}nwall--Bellman inequality.
\begin{lemma}[\cite{ZG17}]\label{lemma_volt}
 {Let vector fields $h_i$ be Lipschitz continuous in a domain $D{\subseteq}\mathbb R^n$, and $h_i\in C^2(D\setminus\Xi;\mathbb R)$, where $\Xi=\{\xi{\in } D{:}h_i(\xi)=0\text{ for }1{\le }i{\le }l\}$, and $ L_{h_j}h_i, L_{h_l}L_{h_j}h_i\in
C(D;\mathbb R^n)$  for all $i,j,l=\overline{1,l}$. If $\xi(t)\in D$, $t\in[0,\tau]$, is a  solution of~$\dot \xi=\sum_{i=1}^l h_i(\xi)w_i(t)$ with $u\in C([0,\tau];\mathbb R^m)$ and $x(0)=x^0\in D$, then  $\xi(t)$ can be represented by the Chen--Fliess series:
\begin{equation}\label{volt1}
 \begin{aligned}
    \xi(t)=\xi^0&+{\sum_{{i_1}=1}^{l}}h_{i_1}(\xi^0)\int\limits_0^t w_{{i_1}}(v)dv+\sum_{\hspace{-0.75em}{i_1},{i_2}=1}^{l} L_{h_{i_2}}h_{i_1}(\xi^0)\int\limits_0^t\int\limits_0^v  w_{{i_1}}(v) w_{{i_2}}(s)dsdv\\
&+R(t),\qquad\qquad\qquad\qquad\qquad t\in[0,\tau],\\
  R(t)=&\sum\limits_{\hspace{-0.75em}{i_1},{i_2},{i_3}=1}^{l}{\int\limits_0^t}{\int\limits_0^v}{\int\limits_0^s} L_{h_{i_3}} L_{h_{i_2}}h_{i_1}(\xi(p))  w_{{i_1}}(v)w_{{i_2}}(s)w_{{i_3}}(p)dpdsdv \end{aligned}
\end{equation}
 is the remainder of the Chen--Fliess series expansion.}
     \end{lemma}
     \begin{lemma}[follows from~\cite{GZE18}]~\label{lemma_rem}
     \\
  { Let the conditions of Lemma~\ref{lemma_x} be satisfied and let $\xi^*\in D$. Assume that there exist $M_1,M_3{\ge}0$, $m{\ge} 1$,  $\varpi\in\{0\}\cup[1,\infty)$ such that
{ $$
\begin{aligned}
\max\limits_{1\le {i_1}\le l}&\|h_{i_1}(\xi(0))\|\le M_1\|\xi(0)-\xi^*\|^m,\\
\max\limits_{1\le{i_1},{i_2},{i_3}\le l}&\|L_{h_{i_3}}L_{h_{i_2}}h_{i_1}(\xi)\|\le M_3\|\xi-\xi^*\|^{\varpi}\text{ for all } \xi\in D.
\end{aligned}
$$}
Then, for all $t\in[0,\tau]$, the remainder $R(t)$ of the Chen--Fliess expansion~\eqref{volt1} of $x(t)$ satisfies the estimate
$$
\begin{aligned}
\|R(t)\|\le &\dfrac{2^{\varpi-2}}{3}(t\nu)^3\|\xi(0)-\xi^*\|^\varpi M_3 \\
&\times\Big(1+M_1(\tau\nu)^\varpi e^{\nu L\varpi\tau}\|\xi(0)-\xi^*\|^{\varpi(m-1)}\Big).
\end{aligned}
$$}
  \end{lemma}
     \begin{lemma}[\cite{GZE18}]~\label{lemma_decay}
   { Let $D\subseteq\mathbb R^n$ be a bounded convex domain,  $W\in C^2(D;\mathbb R)$, $x^*\in D$, and let the following inequalities hold:
   $$
   \begin{aligned}
    & \sigma_{11}\|x-x^*\|^{2m} \le W(x)\le \gamma_{12}\|x-x^* \|^{2m},\\
    &  \sigma_{21} W(x)^{2-\dfrac{1}{m}}\le \|\nabla W(x)\|^2\le\sigma_{22} W(x)^{2-\dfrac{1}{m}},\\
    &\left\|\dfrac{\partial^2 W(x)}{\partial x^2}\right\|\le\sigma_3 W(x)^{1-\dfrac{1}{m}},\\
  %  \|r_\varepsilon\|&\le \Omega\varepsilon^{3/2}V(x^0)^{\dfrac{2m_1(1+m_2)-1}{2m_1}},
    \end{aligned}
    $$
   where $m\ge1$ and  $\sigma_{11},\sigma_{12},\sigma_{21},\sigma_{22},\sigma_3$ are positive constants.
   Then, for any $x^0=x(0)\in D\setminus\{x^*\}$ and any function $x:[0,\varepsilon]\to D$  satisfying the conditions
{    \begin{equation*}%\label{exp_x}
 x(0)=x^0,\;   x(\varepsilon)=x^0-\gamma\varepsilon\nabla W(x^0)+r_\varepsilon,\,\gamma>0,\,r_\varepsilon\in\mathbb R^n,
    \end{equation*}}
the function $W$ satisfies the estimate:
 {\begin{align*}
W(x(\varepsilon))\le W(x^0)\Big(1-\dfrac{\varepsilon\varkappa_1}{m}W^{1-\dfrac{1}{m}}(x^0)+\dfrac{\varepsilon^2\varkappa_2}{2m^2}W^{2-\dfrac{2}{m}}(x^0)\Big)^{m},
 \end{align*}}
where { $\varkappa_1=\gamma\sigma_{21}-{\sqrt{\sigma_{22}}\|r_\varepsilon\|}{ W^{\dfrac{1}{2m}-1}(x^0)}/\varepsilon$,
$\varkappa_2=((m-1)\sigma_{22}+m\sigma_{12})\bigg(\gamma\sqrt{\sigma_{22}}+{\|r_\varepsilon\|}{ W^{\dfrac{1}{2m}-1}(x^0)}/\varepsilon\bigg)^2$.}}
  \end{lemma}
\section{Proof of Theorem~\ref{thm_step1}}    % Each appendix must have a short title.

For the sake of clarity, we divide the proof into several steps resulting in intermediate  statements.
\\
\emph{\underline{Step 0.} Notations and preliminary constructions. }
 To practically stabilize system~\eqref{nonh} at $(0,x^*)$, we will focus on three parameters: $\gamma_1$, $\varepsilon$, and $\mu$, assuming that $\varepsilon<\mu$.   In the proof,  we will determine   big enough $\gamma_1=\gamma_1(\mu)$, small enough $\varepsilon=\varepsilon(\gamma_1,\mu)$, and  small enough $\mu$. It can be seen from the proof that such a choice is always possible.
  %Throughout the proof, we assume that
%\begin{equation}\label{eta}
%\eta\ge\eta_0(\mu,\varepsilon)=\dfrac{1}{\mu}\max\Big\{1,\varepsilon\Big\},
%\end{equation}
We use the following notations in the proof:
for any $\tau\in[0,\varepsilon]$,
{
\begin{align*}
&U(x(\varepsilon),\xi(\varepsilon),\tau)=\max_{\varepsilon\le t\le \varepsilon+\tau}\sum_{i=1}^{m} |\varphi_i^\varepsilon(x(\varepsilon),\xi(\varepsilon),t)|,\\
&  W(\tau)=\max_{0\le t\le \tau}\sum_{j=1}^{2n} |v_j^\tau(t)|\le  \dfrac{c_{w}}{\sqrt\mu},\;c_{w}=2\sum_{j=1}^{n}\sqrt{2\pi k_j}.
\end{align*}}
Recall that the state-dependent control coefficients are defined by~\eqref{a}, which implies that, for any $x(0)=x^0\in D$, $\xi(0)=\xi^0\in D$,
$$
\|a(x^0,\xi^0)\|\le \gamma_1\alpha\|x^0-\xi^0\|.
$$
The H\"{o}lder inequality implies that, for any $\varepsilon>0$ and all $\tau\in[0,\varepsilon]$,
{\begin{equation}\label{est_u}
\begin{aligned}
&U(x^0,\xi^0,\tau)\varepsilon\le \varepsilon \sum_{i\in S_1}|a_i(x^0,\xi^0)|
+2\sqrt{2\pi\varepsilon}\\
&\times\sum_{(i_1,i_2)\in S_2}\sqrt{\kappa_{i_1i_2}|a_{i_1i_2}(x^0,\xi^0)|}\le {c_u}\sqrt {\gamma_1\varepsilon\|x^0-\xi^0\|},
\end{aligned}
\end{equation}
where  $c_{u}= \sqrt{\alpha}(\sqrt{\gamma_1\alpha\varepsilon\|x^0-\xi^0\||S_1|} + 2\sqrt{2\pi})\Big(\sum_{(j_1,j_2)\in S_2}{\kappa^{2/3}_{j_1j_2}}\Big)^{3/4}$} is strictly monotonically increasing w.r.t. $\varepsilon$.

For any {$\delta\in\Big(0,\sqrt{\dfrac{\sigma_{11}}{\sigma_{12}}}{\rm dist}(x^*,\partial D)\Big)$}, let { $\delta_x\in\Big(\sqrt{\dfrac{\sigma_{12}}{\sigma_{11}}}\delta,{\rm dist}(x^*,\partial D)\Big),$} and let  $D'$ be compact,
$D_x=\overline{B_{\delta_x}(x^*)}{\subset} D'{\subseteq} D.$
If $D$ is compact, then we take $D'=D$.
%Let us take any $\delta_x>0$ and $D'\subset \mathbb R^n$ such that
%$$D_x=\overline{B_{\delta_x}(x^*)}\subset D'\subset D.$$
By the conditions of Theorem~1, there exist  $M_f,M_g,M_{3g}>0$  such that, for all $x,\xi\in D_x$,
{
\begin{align}
&\|f_i(x)\|\le M_f,\,\|g_{j}( J(x))\|\le M_g,\,i=\overline{1,m},\ j=\overline{1,2n}\label{M}\\
&\| L_{f_{j_1}}f_{j_2}(x)\|\le M_{2f},\ \big\| L_{f_{j_3}} L_{f_{j_2}}f_{j_1}(x)\big\|\le 6M_{3f},\,j_1,j_2,j_3=\overline{1,m}.\nonumber
%&\|f_i(x)-f_i(\xi)\|\le L_f\|x-\xi\|,\\
%&\|g_{j_1}(J(x))-g_{j_2}(J(\xi))\|\le L_g\|x-\xi\|,\\
%\| L_{g_{j_1}} L_{g_{j_2}}g_{j_3}(J(x))\|\le M_{3g}.
\end{align}
}
If~\eqref{M} and inequalities from the fourth condition of Theorem~1 hold globally in $D$, then we take $D'=D$.
\\
\underline{\emph{Step 1.}} \emph{At this step we construct some a priori estimates which will be exploited further in the proof.}
\\
It is easy to see that  the $\pi_\varepsilon$-solutions of system~\eqref{nonh}  satisfy
{\begin{align}
&\| x(t)- x(0)\|\le M_f{c_u}\sqrt {\gamma_1\varepsilon\|x(0)-\xi(0)\|}\text{ for all }t\in[0,\varepsilon].%,\big(c_{u1}\|x^0-\xi^0\|+c_{u2}\sqrt{\|x^0-\xi^0\|}\big),%c_u\sqrt{\varepsilon\|\bar x(0)-\bar \xi(0)\|},
\label{barx_est}
%\\
%&\| \xi(t)- \xi(0)\|\le \dfrac{\nu\varepsilon}{\eta\sqrt\mu},\label{bary_est}
\end{align}}
Let $\rho_1{>}0$ be  given,
$\nu=c_wM_g$, $\varsigma{>}0$, $\rho_0={\rho_1\mu^\varsigma\sqrt\mu}/{3}$, $
\delta_\xi=\delta_x+\nu\sqrt\mu,$
 $D_\xi=\overline{B_{\delta_\xi}(x^*)}$, $d={\rm dist}(x^*,\partial D')-\delta_x>0, $ and let $\mu_0$ be the smallest positive root of the equation
{\begin{equation}\label{eps0}
  \sqrt\mu\Big({2\rho_1\mu^\varsigma}/{3}+\nu\Big)=d.
\end{equation}}
Obviously, for any $\mu\in(0,\mu_0]$, $\nu\sqrt\mu< d$, so that $\delta_\xi<\delta_x+d$ and $D_\xi\subseteq D'$.
We will also assume that
{\begin{equation}\label{gamma0}
\gamma_1>\bar\gamma_1(\mu)=\dfrac{3\nu}{\rho_1\mu^{\varsigma+1}}.
\end{equation}
}
Such a choice of $\gamma_1$ will be motivated in Step 2.\\
Next, we take
{\begin{equation}\label{vareps0}
\varepsilon_0(\gamma_1,\mu)=\dfrac{1}{\gamma_1}\min\Big\{1,\dfrac{\rho_1\mu^{\varsigma}\sqrt\mu}{3M_f^2c_u^2}\Big\},
\end{equation}}
and observe that
$$
\varepsilon_0(\gamma_1,\mu)\le \dfrac{\rho_1\mu^{\varsigma+1}}{3\nu}
$$
because of~\eqref{gamma0}.

From~\eqref{barx_est} and~\eqref{vareps0} we obtain that,
for each $\mu\in(0,\mu_0]$, $\gamma_1\in(\bar\gamma_1,\infty)$, $\varepsilon\in(0,\varepsilon_0(\gamma_1,\mu)]$,  and for any  $x(0)\in D_x$, $\xi(0)\in \overline{B_{\rho_0}(x(0))}$,  if  $\|\xi(t)-\xi(0)\|\le \dfrac{\nu\varepsilon}{\sqrt\mu}$ with $t\in[0,\varepsilon]$ then
{$$
\begin{aligned}
\| x(t)- x(0)\|&\le M_f{c_u}\sqrt {\gamma_1\varepsilon\|x(0)-\xi(0)\|}\le M_fc_u\sqrt {\gamma_1\varepsilon\rho_0}\\
&\le M_fc_u\sqrt {\dfrac{\gamma_1\varepsilon\rho_1\mu^\varsigma\sqrt\mu}{3}}\le M_fc_u\sqrt {\dfrac{\gamma_1\varepsilon\rho_1\mu^\varsigma\sqrt\mu}{3}}\le\dfrac{\rho_1\mu^\varsigma\sqrt\mu}{3},\\
\|x(t)-\xi(t)\|&\le \|x(t)-x^0\|+\|x^0-\xi^0\|+\|\xi(t)-\xi^0\| \\
&\le M_f{c_u}\sqrt {\dfrac{\gamma_1\varepsilon\rho_1\mu^\varsigma\sqrt\mu}{3}}+\dfrac{\rho_1\mu^\varsigma\sqrt\mu}{3}+\dfrac{\nu\varepsilon}{\sqrt\mu}\le \rho_1\mu^\varsigma\sqrt\mu.
\end{aligned}
$$}
If, additionally, $\|\xi^0-\xi^*\|\in D_\xi$ then
{$$
\begin{aligned}
\|x(t)-x^*\|&\le \|x(t)-x^0\|+\|x^0-\xi^0\|+\|\xi^0-\xi^*\|\\
&\le \dfrac{2\rho_1\mu^\varsigma\sqrt\mu}{3}+\delta_\xi\le {\rm dist}(x^*,\partial D').
%\| \xi(t)- \xi(0)\|&\le \nu\sqrt\mu \le e{\rm dist}(x^*,\partial D')\text{ for }t\in[0,\eta\mu],
\end{aligned}
$$}
 This proves the following intermediate statement.
{\begin{statement}
 { For any $\mu\in(0,\mu_0]$, $\gamma_1\in(\bar\gamma_1,\infty)$, $\varepsilon\in(0,\varepsilon_0(\gamma_1,\mu)]$, $x^0\in D_x$, the $\pi_\varepsilon$-solutions of system~\eqref{nonh} with the initial conditions $x(0)=x^0,\xi(0)=\xi^0$ satisfy the following property:
  $$
\|x^0-\xi^0\|\le\dfrac{\rho_1\mu^\varsigma\sqrt\mu}{3}{\Rightarrow}\|x(t)-\xi(t)\|\le\rho_1\mu^\varsigma\sqrt\mu\text{ for }t\in[0,\varepsilon].
  $$
Furthermore, if  $\xi^0\in D_\xi\subseteq D'$ then $x(t)$ is well-defined in $D'$ for  $t\in[0,\mu]$.
} \end{statement}}

\underline{\emph{Step 2.}}
\emph{Our next goal is to ensure that the $x$-component of the $\pi_\varepsilon$ solution of system~\eqref{nonh} is in a sufficiently small neighborhood of the $\xi$-component.}

For this, we apply Lemma~\ref{volt1}. Namely, assume  that $x(t)\in D_x$ for $t\in[0,\varepsilon]$, $\xi(0)=\xi^0\in B_{\rho_1}(x^0)$. Then
{\begin{align}
\| \xi(t)- \xi(0)\|\le \dfrac{\nu\varepsilon}{\sqrt\mu},\label{bary_est}
\end{align}}
and
the $\pi_\varepsilon$-solution $ x(t)$ of system~\eqref{nonh} with controls~\eqref{cont}  can be represented my means of  the  Chen--Fliess series:
{\begin{equation}\label{xeps}
     \begin{aligned}
 x(\varepsilon)= x^0-  \varepsilon \gamma_1(x^0-\xi^0)+R_1(\varepsilon)+ R_2(x^0,\xi^0,\varepsilon),
\end{aligned}
\end{equation}}
where $R_1(\varepsilon)$ is defined from Lemma~\ref{lemma_volt}, and
{\begin{align*}
 R_2(x^0,&\xi^0,\varepsilon)={\varepsilon^{3/2}}\sum_{j_1\in S_1}\sum_{j_2=1}^m[f_{j_1},f_{j_2}](x^0) a_{j_1}(x^0,\xi^0)\sum_{q:(q,j_2)\in S_2}
 \sqrt{\dfrac{|a_{qj_2}(x^0,\xi^0)|}{\pi K_{qj_2}}}
 \\
 &+\dfrac{\varepsilon^2}{2}\sum_{j_1,j_2\in S_1} L_{f_{j_2}}f_{j_1}(x^0) a_{j_1}(x^0,\xi^0)a_{j_2}(x^0,\xi^0).
  \end{align*}}
  Denote $R(x^0,\xi^0,\varepsilon)=R_1(\varepsilon)+R_2(x^0,\xi^0,\varepsilon)$.
Using~\eqref{est_u} and notations from~\eqref{M}, we get
$$
     \|R_1(\varepsilon)\|\le M_{3f}c_u^3\big(\varepsilon\|x^0-\xi^0\|\big)^{3/2}\text{ for all }t\in[0,\varepsilon],
$$
and
  \begin{align}
  \| R(x^0,\xi^0,\varepsilon)\|\le \zeta_1\big(\varepsilon\|x^0-\xi^0\|\big)^{3/2},\label{est_R} %\le \sigma_1 \varepsilon^{3/2}\|x(\tau_0)-\xi(\tau_0)\|^{3/2},
\end{align}
$$
\zeta_1= M_{3f}c_u^3+\dfrac{ M_{2f}}{2}\sqrt{\nu\varsigma\alpha}(\gamma_1\alpha)^{3/2} +2(\gamma_1\alpha)^{3/2}M_{2f} \sqrt{|S_1|}\sum_{j_1=1}^m\Big(\sum_{(j_2,j_1)\in S_2}\kappa_{j_2j_1}^{-2/3}\Big)^{3/4}.
$$
Combining~\eqref{est_R}, \eqref{bary_est}, and~\eqref{xeps}, we come to the following estimate:
%First, we estimate
%$$
%\|x(\varepsilon)-x^0\|\le \gamma_1\varepsilon+\zeta_1\big(\varepsilon\|x^0-\xi^0\|\big)^{3/2}\le  \gamma_1\varepsilon+\zeta_1\big(\varepsilon\rho_1\big)^{3/2}.
%$$
%Taking $\varepsilon_1$ as the smallest positive root of the equation
%$$
%\gamma_1\varepsilon+ \zeta_1\big(\varepsilon\rho_1\sqrt\mu\big)^{3/2}=\delta_x,
%$$
%we ensure that if $\|x^0-\xi^0\|\le \rho_1\sqrt\mu$, then $\|x(\varepsilon)-x^0\|\le \delta_x$ for any $\varepsilon\in(0,\min\{\varepsilon_0,\varepsilon_1\})$, i.e. $x(\varepsilon)\in D_x$.
%
{$$
\|x(\varepsilon)-\xi(\varepsilon)\|\le(1-\varepsilon\gamma_1)\|x^0-\xi^0\|+\zeta_1\big(\varepsilon\|x^0-\xi^0\|\big)^{3/2}+\dfrac{\nu\varepsilon}{\sqrt\mu}.
$$}
For any $\gamma_1>\bar\gamma_1$, let $\lambda_1\in[\bar\gamma_1,\gamma_1)$ and define
{\begin{equation}\label{vareps1}
\varepsilon_1(\gamma_1,\mu)=\min\Big\{\varepsilon_0(\mu),\Big(\dfrac{\gamma_1-\lambda_1}{\zeta_1\sqrt\delta_x}\Big)^2\Big\}.
\end{equation}}
 Recall that $\varepsilon\lambda_1<\varepsilon\gamma_1<1$.
 Then, for any $\varepsilon\in(0,\varepsilon_1(\gamma_1,\mu_1))$,
{$$
\begin{aligned}
\|x(\varepsilon)-\xi(\varepsilon)\|&< (1-\varepsilon\bar\gamma_1)\|x^0-\xi^0\|+\dfrac{\nu\varepsilon}{\sqrt\mu}.
\end{aligned}
$$}
Recall that $\bar\gamma_1$ is given by~\eqref{gamma0}, which implies
$
\dfrac{\nu\varepsilon}{\sqrt\mu}=\dfrac{\gamma_1\rho_1\mu^\varsigma\sqrt\mu}{3}.
$
 This together with Statement~1 gives us the next intermediate result.
{\begin{statement}
  {Assume that $x(t)\in D'$ for all $t\in[0,\varepsilon_0]$, $x(0)\in D_x$. Then,  for any $\mu\in(0,\mu_0]$, $\gamma_1\in(\bar\gamma_1,\infty)$, $\varepsilon\in(0,\varepsilon_1(\gamma_1,\mu)]$, the following properties hold:
$$
\begin{aligned}
%\text{if }\|x^0-\xi^0\|\ge \dfrac{\rho_1\sqrt\mu}{2}\text{ then }\|x(\varepsilon)-\xi(\varepsilon)\|&<\|x^0-\xi^0\|;\\
\text{if }&\|x^0-\xi^0\|\le \dfrac{\rho_1\mu^\varsigma\sqrt\mu}{3}\text{ then }\|x(\varepsilon)-\xi(\varepsilon)\|\le\dfrac{\rho_1\mu^\varsigma\sqrt\mu}{3},\\
\text{ and }&\|x(t)-\xi(t)\|\le \rho_1\mu^\varsigma\sqrt\mu\text{ for all }t\in[0,2\varepsilon].
\end{aligned}
$$}
\end{statement}
}
\underline{\emph{Step 3.}}
Now let us put $x(0)=x^0=\xi^0=\xi(0)$, $x^0\in D_x$. Then $x(t)\equiv x^0\in D_x$ for all $t\in[0,\varepsilon]$, and
{$$
\|\xi(t)-\xi^0\|\le \|\xi(t)-\xi^0\|+\|\xi^0-\xi^*\|\le\nu \sqrt\mu +\delta_x=\delta_\xi,
$$}
i.e. $\xi(t)\in D_\xi\subset D'$  for $t\in[0,\varepsilon]$. Besides, Statement~2 implies
$$
\|x(\varepsilon)-\xi(\varepsilon)\|<\dfrac{\rho_1\mu^\varsigma\sqrt\mu}{3}.
$$
From Statements~1 and~2, the $x$-component of the $\pi_\varepsilon$-solution of system~\eqref{nonh} is also well-defined in $D'$ for $t\in[\varepsilon,2\varepsilon]$. Again, it is easy to see that
$$
\|\xi(t)-\xi^0\|\le \nu \sqrt\mu +\delta_x\text{ for }t\in[0,2\varepsilon],
$$
 i.e.  $\xi(2\varepsilon)\in D_\xi$ and
 $$
\|x(2\varepsilon)-\xi(2\varepsilon)\|<\dfrac{\rho_1\mu^\varsigma\sqrt\mu}{3}.
$$
Without loss of generality, we may assume that  $\dfrac{\mu}{\varepsilon}=\mathcal N_1,\text{ with some }\mathcal N_1\in\mathbb N.$
Repeating Steps~1--2 until $t=\mathcal N\varepsilon$, we come to the following statement.
\begin{statement}
 For any  $\mu\in(0,\mu_0]$, $\gamma_1\in(\bar\gamma_1,\infty)$, $\varepsilon\in(0,\varepsilon_1(\gamma_1,\mu)]$, the $\pi_\varepsilon$-solutions $(x(t),\xi(t))$ of  system~\eqref{nonh} with the initial conditions $x(0)=\xi(0)\in D_x$ are well-defined in $D'\times D'$ for all $t\in[0,(\mathcal N_1+1)\varepsilon]$,
$$
\|x(t)-\xi(t)\|\le \rho_1\mu^\varsigma\sqrt\mu\text{ for all }t\in[0,\mu],\ \|x(\mu)-\xi(\mu)\|\le \dfrac{\rho_1\mu^\varsigma\sqrt\mu}{3}.
$$
\end{statement}
Thus, for any  $\mu\in(0,\mu_0]$, we can take $\gamma_1(\mu)$, $\varepsilon(\gamma_1(\mu),\mu)$, such that   $x(t),\xi(t)\in D'$ for $t\in[0,\mu]$.
In the next steps, we will find sufficiently small $\mu$ independently on $\varepsilon$ and $\gamma_1$.

\underline{\emph{Step 4.}} \emph{The goal of this step is to ensure the decay of the cost function $J(x)$ along the trajectories of system~\eqref{nonh} by choosing sufficiently small $\mu$.}

For this purpose we apply again Lemma~\ref{volt1}.
Since $x(t),\xi(t)\in D'$ for $t\in[0,\mu]$, we may consider the Chen--Fliess series expansion of the $\xi$-component of solution of system~\eqref{nonh} on the interval $[0,\mu]$:
{ \begin{align}\label{volt_y}
\xi(\mu)&=\xi^0-\mu\gamma_2\nabla J(\xi^0)+R_3(\mu),
\end{align}}
where
{\begin{align*}
R_3(\mu)=&\sum_{j=1}^{2n}\int_0^{\mu}\Big(g_{j}( J(x(s_1)))-g_{j}( J(\xi(s_1)))\Big)e_jv_{j}^\mu({s_1})ds_1\\
 &+\sum_{j_1,j_2=1}^{2n}\int_0^{\mu}\int_0^{s_1} L_{e_{j_2}\big(g_{j_2}\circ J(x(s_2))-g_{j_2}( J(\xi(s_2)))\big)} g_{j_1}( J(\xi(s_2)))e_{j_1}v_{j_2}^\mu ({s_2})v_{j_1}^\mu({s_1})ds_2ds_1\\
& +\sum_{j_1,j_2,j_3=1}^{2n}\int_0^{\mu}\int_0^{s_1}\int_0^{s_2} v_{j_3}^\mu ({s_3})v_{j_3}^\mu ({s_2})v_{j_1}^\mu({s_1}) \\
&\times L_{e_{j_3}g_{j_3}( J(x(s_3)))} L_{e_{j_2}g_{j_2}( J(\xi(s_3)))}g_{j_1}( J(\xi(s_3)))e_{j_1}ds_3ds_2ds_1.
  \end{align*}}
  Under the assumptions of Theorem~1, we conclude that
 { $$
\begin{aligned}
  \|R_3(\mu)\|\le c_w{\sqrt\mu}(L_g+\sqrt\mu L_{2g}c_w)\max\limits_{0\le t\le\mu}\|x(t)&-\xi(t)\|+\mu^{3/2}M_{3g}.
\end{aligned}
  $$}
Thus, applying Statement~3 we get
$$
\|R_3(\mu)\| \le \zeta_2 \mu^{1+\tilde\varsigma},
$$
where $\tilde\varsigma=\min\{\varsigma,1/2\}$, $\zeta_2=c_w\mu^{\max\{0,\varsigma-1/2\}}\rho_1(L_g+\sqrt\mu L_{2g}c_w)+\mu^{\max\{0,1/2-\varsigma\}}M_{3g}$.
\\
Using Taylor's formula for the function $J(\xi)$,
{$$
\begin{aligned}
J&(\xi(t))=J(\xi^0)+\nabla J(\xi^0)(\xi(t)-\xi^0)+\dfrac{1}{2}\sum_{i,j=1}^{2n}\dfrac{\partial^2 J(x)}{\partial x_i\partial x_j}\Big|_{x=\xi^0+\theta \xi(t)}(\xi_i-\xi_i^0)(\xi_j-\xi_j^0),
\end{aligned}
$$}
and exploiting~\eqref{J}, we obtain
{$$
\begin{aligned}
J(\xi(\mu))&\le J(\xi^0)-\mu\gamma_2\sigma_{21}J(\xi^0)+\mu^{1+\varsigma}\zeta_2\sqrt{\sigma_{22}J(\xi^0)}+\sigma_3\big(\mu^2\gamma_2^2\sigma_{22}^2J(\xi^0)+\zeta_2^2\mu^{2+2\varsigma}\big)\\
&=J(\xi^0)\Big(1-\mu\gamma_2\big(\sigma_{21}-\mu\gamma_2\sigma_3\sigma_{22}^2\big)\Big)+\mu^{1+\varsigma}\zeta_2\Big(\sqrt{\sigma_{22}J(\xi^0)}+\sigma_3\zeta_2\mu^{1+\varsigma}\Big).
\end{aligned}
$$}
Let $\mathcal L_c=\{x\in D:J(x)\le c\}$, $c_J=\sigma_{11}\delta_x^2$. Then
{\begin{equation}\label{lset}
  \overline{B_\delta(x^*)}\subseteq \mathcal L_{c_J} \subseteq D_x.
\end{equation}}
For any $\rho_2\in(0,c_J]$, $\lambda_2\in(0,\gamma_2\sigma_{21})$, we define
{\begin{equation}\label{eps1}
\mu_1=\min\{\mu_0,{1}/{\lambda_2},\hat\mu_1\},
\end{equation}}
where $\hat\mu_1$ is the smallest positive root of the equation
{$$
\rho_2\mu\gamma_2\sigma_3\sigma_{22}^2+{\mu^{\varsigma}}\zeta_2\Big(\sqrt{\sigma_{22}\rho_2}+\sigma_3\zeta_2\mu^{1+\varsigma}\Big)=\rho_2(\gamma_2\sigma_{21}-\lambda_2).
$$}
Then, for any $\mu\in(0,\mu_1)$, the following two scenarios are possible:

S1) If { $J(\xi^0)\le\rho_2$} then {$J(\xi(\mu))\le \rho_2(1-\mu\lambda_2){<}\rho_2$.}
In this case, $\xi(\mu)\in D_x$, Additionally, Statement~3 implies that  $\|x(\mu)-\xi(\mu)\|\le \dfrac{\rho_1\mu^\varsigma\sqrt\mu}{3}$. Repeating the above argumentation, we get $\xi(N\mu)\in D_x$ for all natural numbers $N$.

S2) If $J(\xi^0)>\rho_2$ then $J(\xi(\mu))<J(\xi^0)(1-\mu\lambda_2)<J(\xi^0)$.

Consider S2).
If $\xi^0=x^0\in \overline{B_\delta(x^*)}$ then $\xi(\mu)\in  \mathcal L_{c_J} \subseteq D_x$. Again, Statement~3 gives  $\|x(\mu)-\xi(\mu)\|\le \dfrac{\rho_1\mu^\varsigma\sqrt\mu}{3}$. Thus, we may repeat all the steps for $t\in[\mu,2\mu]$.

Summarizing all the above, we arrive at the following conclusion:
there exists an $\mathcal N_2\in\mathbb N\cup\{0\}$ such that
{$$
\begin{aligned}
&J(\xi(t))\le J(\xi^0)e^{-\lambda_2t}\text{ for }t=0,\mu,\dots, (\mathcal N_2-1)\mu,\\
&J(\xi(t))\le \rho_2\text{ for }t=\mathcal N_2\mu,(\mathcal N_2+1)\mu,\dots  \, .
\end{aligned}
$$}
Consequently,  
$$\|\xi(t)-x^*\|\le \sqrt{\dfrac{\sigma_{12}}{\sigma_{11}}}\|x^0-x^*\|e^{-\lambda_2t}\text{ for }t=0,\mu,\dots,(\mathcal N_2-1)\mu,$$
$$\|\xi(\mathcal N_2\mu)-x^*\|\le \sqrt{\dfrac{\rho_2}{\sigma_{11}}}\le \delta_x\text{ for  }t=\mathcal N_2\mu,(\mathcal N_2+1)\mu,\dots .$$
For an arbitrary $t\in[0,\mathcal N_2\mu]$, we denote the integer part of $\dfrac{t}{\mu}$ as $\Big[\dfrac{t}{\mu}\Big]$ and observe that $0<t-\Big[\dfrac{t}{\mu}\Big]\mu<\mu$. Then
{$$
\begin{aligned}
\|\xi(t)-x^*\|&\le \Big\|\xi\Big(\Big[\dfrac{t}{\mu}\Big]\mu\Big)-x^*\Big\|+\Big\|\xi(t)-\xi\Big[\dfrac{t}{\mu}\Big]\Big\|\\
&\le \sqrt{\dfrac{\sigma_{12}}{\sigma_{11}}}\|x^0-x^*\|e^{-\lambda_2\Big[\dfrac{t}{\mu}\Big]\mu}+\nu\sqrt\mu\le \beta\|x^0-x^*\|e^{-\lambda_2 t}+\nu\sqrt\mu,
\end{aligned}
$$}
where $\beta=\sqrt{\dfrac{\sigma_{12}}{\sigma_{11}}}e^{\lambda_2\mu}$. This yields the following result.
{\begin{statement}
 {For any $\mu\in(0,\mu_1]$, $\gamma_1\in(\bar\gamma_1,\infty)$, $\varepsilon\in(0,\varepsilon_0(\gamma_1,\mu)]$,  the $\pi_\varepsilon$-solutions $(x(t),\xi(t))$ of  system~\eqref{nonh} with the initial conditions $x(0)=\xi(0)\in D_x$ are well-defined in $D'\times D'$ for all $t\in[0,\infty)$, and the following estimates hold:
$$
\begin{aligned}
&\|\xi(t)-x^*\|\le \beta\|x^0-x^*\|e^{-\lambda t}+\nu\sqrt\mu\text{ for }t\in[0,\mathcal N_2\mu],\\
&\|\xi(t)-x^*\|\le \sqrt{\dfrac{\rho_2}{\sigma_{11}}}+\nu\sqrt\mu\text{ for }t\in[\mathcal N_2\mu,\infty).
\end{aligned}
$$
Furthermore,
$$
\|x(t)-\xi(t)\|\le \rho_1\mu^\varsigma\sqrt\mu\text{ for all }t\in[0,\infty).
$$}
\end{statement}}

\underline{\emph{Step 5.}}
 \emph{Finally, we estimate $\|x(t)-x^*\|$ for $t\in[0,\infty)$}.
 
Applying the triangle inequality together with Statement~4, we get the following:
{$$
\begin{aligned}
&\|x(t)-x^*\|\le \beta\|x^0-x^*\|e^{-\lambda t}+\rho_1\mu^\varsigma\sqrt\mu+\nu\sqrt\mu\text{ for }t\in[0,\mathcal N_2\mu],\\
&\|x(t)-x^*\|\le \rho_1\mu^\varsigma\sqrt\mu+\sqrt{\dfrac{\rho_2}{\sigma_{11}}}+\nu\sqrt\mu\text{ for }t\in[\mathcal N_2\mu,\infty).
\end{aligned}
$$}
Since $\rho_1$, $\rho_2$ are arbitrary and $\mu$ can be chosen small enough,
the above inequalities imply the assertion of Theorem~1.
In particular, for an arbitrary $\rho>0$, one can take
$\rho_1>0$ and $\mu>0$ such that
\begin{equation}\label{eps2}
\rho_1\mu^\varsigma\sqrt\mu+\nu\sqrt\mu\le\dfrac{\rho}{2},
\end{equation}
 and $\rho_2\le\dfrac{1}{4}\rho^2\sigma_{11}$. Then
$$
\|x(t)-x^*\|\le \beta\|x^0-x^*\|e^{-\lambda t}+\rho\text{ for all }t\in[0,\infty).
$$
Note that the choice of $\mu$ does not depend on $\varepsilon,\eta$, and the choice of $\gamma_1$ does not depend on $\varepsilon$. Namely, given $\delta,\rho,\rho_1,\rho_2$, one can choose a $\bar\mu>0$ satisfying~\eqref{eps0},~\eqref{eps1} and~\eqref{eps2},   and take any $\hat\mu\in(0,\bar\mu]$. The next step is to determine $\bar\gamma_1(\hat\mu)$
 satisfying~\eqref{gamma0}, and take any $\hat\gamma_1\in(\bar\gamma_1,\infty)$. Finally, $\bar\varepsilon(\hat\gamma_1,\hat\mu)$ has to be specified according to~\eqref{vareps0} and~\eqref{vareps1}.
\end{document}